\documentclass[12pt,twoside]{amsart}

\usepackage{amsmath, amsthm, amssymb, graphicx, paralist, bm}

\evensidemargin \oddsidemargin
\newtheorem{thm}{Theorem}[section]

\newtheorem{cor}[thm]{Corollary}

\newtheorem{prop}[thm]{Proposition}

\theoremstyle{definition}

\newtheorem{definition}[thm]{Definition}

\theoremstyle{remark}
\newtheorem{remark}[thm]{Remark}

\newtheorem{example}[thm]{Example}

\numberwithin{equation}{section}

\title[Split Orders and Convex Polytopes in Buildings] {Split Orders and Convex Polytopes in Buildings} 

\author {Thomas R. Shemanske}

\address{Department of Mathematics\\ 
6188 Kemeny Hall\\
Dartmouth College\\
Hanover, NH 03755}

\email[] {thomas.r.shemanske@dartmouth.edu}
\urladdr{http://www.math.dartmouth.edu/~trs/ }

\date{August 26, 2008}

\newcommand{\A}{{\mathcal A}}
\newcommand{\C}{{\mathcal C}}
\renewcommand{\O}{{\mathcal O}}

\newcommand{\bmu}{{\bm \mu}}
\newcommand{\be}{{\bm e}}
\newcommand{\bb}{{\bm b}}

\newcommand{\bnu}{{\bm \nu}}
\newcommand{\p}{{\mathfrak p}}
\newcommand{\R}{\mathbb{R}}
\newcommand{\Z}{\mathbb{Z}}
\newcommand{\LO}{{\Lambda_0}}

\newcommand{\ds}{\displaystyle}

\renewcommand{\slash}{\,|\,}
\newcommand{\diag}{\mathrm {diag}}
\newcommand{\End}{\mathrm {End}}

\newcommand{\spmat}[4]{\left(\begin{smallmatrix}{#1}&{#2}\\{#3}&{#4}\end{smallmatrix}\right)}

\newcommand{\pmatt}[9]{\begin{pmatrix}{#1}&{#2}&{#3}\\{#4}&{#5}&{#6}\\{#7}&{#8}&{#9}\end{pmatrix}}
\newcommand{\spmatt}[9]{\left(\begin{smallmatrix}{#1}&{#2}&{#3}\\{#4}&{#5}&{#6}\\{#7}&{#8}&{#9}\end{smallmatrix}\right)}

\parskip=\medskipamount

\begin{document}



\maketitle

\begin{abstract}
  As part of his work to develop an explicit trace formula for Hecke
  operators on congruence subgroups of $SL_2(\Z)$, Hijikata
  \cite{Hijikata-traces} defines and characterizes the notion of a
  split order in $M_2(k)$, where $k$ is a local field.  In this paper,
  we generalize the notion of a split order to $M_n(k)$ for $n>2$ and
  give a natural geometric characterization in terms of the affine
  building for $SL_n(k)$.  In particular, we show that there is a
  one-to-one correspondence between split orders in $M_n(k)$ and a
  collection of convex polytopes in apartments of the building such
  that the split order is the intersection of all the maximal orders
  representing the vertices in the polytope.  This generalizes the
  geometric interpretation in the $n=2$ case in which split orders
  correspond to geodesics in the tree for $SL_2(k)$ with the split
  order given as the intersection of the endpoints of the geodesic.
\end{abstract}


\section{Introduction}\label{sec:introduction}

The study of orders in noncommutative algebras has a long history with
known applications to class field theory, modular forms, and geometry.
In \cite{Hijikata-traces}, Hijikata defines and characterizes split
orders in $M_2(k)$, $k$ a local field, as part of his work to develop
an explicit trace formula for Hecke operators on congruence subgroups
of $SL_2(\Z)$. His characterization of split orders is entirely algebraic, 
characterizing them as either maximal orders or the
intersection of two uniquely determined maximal orders.  More
precisely, he shows that

\begin{prop}Let $k$ be a local field, $\O$ its valuation ring, and
  $\p$ the unique maximal ideal of $\O$.  Let $S$ be an $\O$-order in
  $A = M_2(k)$; the following are equivalent and define the notion of
  a {\em split order} in $A$.
  \begin{compactenum}
  \item $S$ contains a subset which is $A^\times$-conjugate to
    $\spmat{\O}{0}{0}{\O}$.
  \item $S$ is $A^\times$-conjugate to $\spmat{\O}{\O}{\p^\nu}{\O}$
    for some non-negative integer $\nu$.
  \item $S$ is the intersection of at most two maximal orders in $A$.
  \item $S$ is either maximal or the intersection of two uniquely
    determined distinct maximal orders.
  \end{compactenum}
\end{prop}

Hijikata's proposition has the following geometric interpretation.
For $k$ a local field, the vertices of the affine building associated
to $SL_2(k)$ are in one-to-one correspondence with the maximal orders
in $M_2(k)$.  Moreover, it is well-known that $SL_2$-building is
actually a $(q+1)$-regular tree ($q$ the cardinality of the residue
field of $k$), so that any two vertices determines a unique path or
geodesic between them.  Thus split orders are in one-to-one
correspondence with the geodesics of finite (nonnegative) length on
the tree, with the split order realized as the intersection of the
maximal orders representing the endpoints of the geodesic; geodesics
of length zero are the maximal orders.

In this paper we consider the generalization of the notion of a split
order to $B = M_n(k)$ for $n > 2$, and give a geometric
characterization using the affine building for $SL_n(k)$.  We take as
a definition of a split order, an order in $M_n(k)$ which contains a
subring $B^\times$-conjugate to $R = \begin{pmatrix}
  \O&&0\\
  &\ddots&\\
  0&&\O\\
\end{pmatrix}$, hereafter denoted as $R= \diag(\O, \dots, \O)$.  The
geometric generalization which we derive agrees with the $n=2$
characterization, though in a manner slightly more nuanced than the
one given above.  In generalizing, one finds that there is no
comparable uniqueness statement (as in Hijikata's proposition) which
characterizes a split order as the intersection of a uniquely
determined minimal set of maximal orders. Rather, the uniqueness
arises by considering the set of all maximal orders which contain the
split order.  We show that
\begin{compactitem}
\item there is an apartment which contains the set of all maximal
  orders containing a given split order,
\item this collection of maximal orders consists of the set of
  all vertices which lie in a convex polytope uniquely
determined by the split order, and
\item the split order is the intersection of all the maximal
  orders in this convex polytope.
\end{compactitem}
In the case $n=2$ (where the building is a tree), Hijikata's result
shows that a split order is the intersection of the maximal orders
which are the endpoints of the geodesic which characterize it; from
this work it follows that the split order is also the intersection of
all the maximal orders contained in the geodesic.  For $n > 2$ and in
the case where all the maximal orders are vertices of a single chamber
in the building, the notion of split orders reduces to that of chain
orders studied (to a different end) in \cite{Abramenko-Nebe}. In that
case the notion of convexity is implicit in the structure of the
building as the convex polytopes which arise are simply faces of the
chamber.  The present work addresses finitely many maximal orders
chosen arbitrarily in any apartment of the building.  The fact that
there is a given apartment containing all the maximal orders which
contain a split order is an interesting extension of the standard
building fact that any two simplicies in a building are contained in a
single apartment, and may point to more complicated structure implicit
in the building.

\textbf{Acknowledgments:} Finally, the author thanks his colleagues
Vladimir Chernov, Scott Pauls, David Webb, and Peter Winkler for
useful discussions concerning combinatorial and geometric aspects of
this project.

\section{Split Orders}
\subsection{Definition and initial characterization}

Let $k$ be a local field, $\O$ its valuation ring, and $\p = \pi\O$
the unique maximal ideal of $\O$, with $\pi$ a fixed uniformizing
parameter.  Let $B$ be the central simple algebra $M_n(k)$, and fix a
subring $R$ having the form $R= \diag(\O, \dots, \O)$.  Recall that an
order $S \subset B$ is a subring of $B$ containing the identity which
is also a free $\O$-module having rank $n^2$.  We begin our
investigation of split orders with the special case in which the order
$S \subset B$ actually contains the subring $R$. We shall see that the
consideration of general split orders (containing a conjugate of $R$)
simply amounts to a change of basis and shifts the geometric
perspective from one apartment to another.

We first give an initial, though somewhat unsatisfying, algebraic
characterization of these split orders. Let $E^{(i,j)}$ denote the $n\times
n$ matrix with a 1 in the $(i,j)$ position and zeros elsewhere.

\begin{prop}\label{prop:entries-form-ideals}
  \begin{compactenum}
  \item Let $S \subset M_n(k)$ be a ring containing $E^{(i,i)}$ for $1
    \le i \le n$.  Then $A = (a_{ij}) \in S$ if and only if
    $a_{ij}E^{(i,j)} \in S$ for all $i,j$.
  \item Let $S$ be an order in $M_n(k)$ containing $E^{(i,i)}$ for $1
    \le i \le n$.  Then $S$ has the form $S =
    \left(\begin{smallmatrix}
        \O&&\p^{\nu_{ij}}\\
        &\ddots&\\
        \p^{\nu_{ij}}&&\O\\
      \end{smallmatrix}\right)$ which we simplify to $S =
    (\p^{\nu_{ij}})$ with the understanding that $\nu_{ii} = 0$ for
    all $i$.
  \item Let $S = (\p^{\nu_{ij}})\subset M_n(k)$ be a set with
    $\nu_{ii} = 0$ for all $i$.  Then $S$ is an order if and only if
    $\nu_{ik} + \nu_{kj} \ge \nu_{ij}$ for every $i, j, k$.
  \end{compactenum}
 
\end{prop}

\begin{proof}
  For the first item, one direction is obvious and for the other,
  simply observe that $E^{(i,i)} A E^{(j,j)} = a_{ij}E^{(i,j)}$.  For
  (2), let $S_{ij} = \{E^{(i,i)} A E^{(j,j)} = a_{ij}E^{(i,j)} \mid A
  \in S\}$.  Since $S$ is an order and hence has rank $n^2$ as an
  $\O$-module, it follows that $S_{ij} \ne \{0\}$ .  Since $S$
  contains all the $E^{(i,i)}$, it is obvious that $S_{ij}$ is a
  fractional $\O$-ideal, hence has the form $\p^{\nu_{ij}}E^{(i,j)}$.
  Since $\O E^{(i,i)} \subseteq S_{ii}$, it is easy to deduce (e.g.,
  from the integrality of elements of $S$ \cite{Reiner-book}) that
  $S_{ii} = \O E^{(i,i)}$.  For (3), if $S$ is closed under
  multiplication, then $S_{ik} S_{kj} \subseteq S_{ij}$, hence
  $\p^{\nu_{ik}} \p^{\nu_{kj}} \subseteq \p^{\nu_{ij}}$, so
  $\nu_{ik} + \nu_{kj} \ge \nu_{ij}$.  Conversely a set $S =
  (\p^{\nu_{ij}})$ with $\nu_{ii} = 0$ is an order if and only if it
  is closed under multiplication.  Let $A = \sum_{i,j}
  a_{ij}E^{(i,j)}$, $B = \sum_{k,\ell} b_{k\ell}E^{(k,\ell)} \in
  S$. Now $AB = \sum_{i,j,k,\ell} a_{ij}b_{k\ell}E^{(i,j)}E^{(k,\ell)}
  = \sum_{i,j,\ell} a_{ij}b_{j\ell}E^{(i,j)}E^{(j,\ell)} =
  \sum_{i,\ell}\left(\sum_j a_{ij}b_{j\ell}\right)E^{(i,\ell)} $.
  Since $a_{ij} \in \p^{\nu_{ij}}$, and $b_{j\ell}\in
  \p^{\nu_{j\ell}}$, the condition $\nu_{ij} + \nu_{j\ell} \ge
  \nu_{i\ell}$ shows that $\sum_j a_{ij}b_{j\ell} \in
  \p^{\nu_{i\ell}}$, and hence $AB \in S$.
\end{proof}

\subsection{The maximal orders which contain a split order}

Next we consider the extent to which the alternate characterizations
of split orders in $M_2(k)$ given by Hijikata hold in $B=M_n(k)$ when
$n > 2$. Naive conjectures concerning a minimal set of maximal orders
whose intersection produces the split order are easily shown not to
hold in general, however a uniqueness statement can be deduced
characterizing split orders as the intersection of a geometrically
distinguished collection of maximal orders which nicely generalizes
the situation for $n=2$.

In particular, we consider whether a split order is characterized by
the set of all maximal orders which contain it.  To that end, we let
$\Lambda_0 = M_n(\O)$ be a fixed maximal order in $B$.  It is well
known \cite{Reiner-book} that every maximal order in $B$ is conjugate
by an element of $B^\times$ to $\Lambda_0$.

We first characterize those maximal orders which contain
the subring $R$, which reduces to characterizing those $\xi = B^\times$, so that
$R \subset \xi^{-1} \LO \xi$.  Since $M_n(k) = k^\times M_n(\O)$ and
the action by conjugation of $k^\times$ is trivial, we may assume that $\xi \in
M_n(\O)$, and in particular, we may choose for $\xi$ any
representative of $GL_n(\O)\xi$.  Thus there is no loss of generality
to assume that $\xi$ is in Hermite normal form (see e.g., \cite{Newman}),
that is $\xi =
\begin{pmatrix}
  \pi^{m_1}&a_{12}&& \dots && a_{1n}\\
  0&\pi^{m_2}&a_{23}& \dots && a_{2n}\\
  0&0&\ddots&&\\
  0&0&\dots&&\pi^{m_{n-1}}&a_{n-1\; n}\\
  0&0&&\dots&&\pi^{m_n}
\end{pmatrix}
$, an upper triangular matrix with powers of the fixed uniformizer on
the diagonal and entries $a_{ij}$ ($i<j$) in a fixed set of residues
of $\O/\pi^{m_j}\O$. We may and do assume the representative of the
zero class is actually zero.

\begin{prop}\label{prop:HNF-to-diagonal}With the notation and assumptions as above, we
  have the $R \subset \xi^{-1} \LO \xi$ if and only if $\xi$ is
  diagonal, $\xi = \diag(\pi^{m_1}, \dots, \pi^{m_n})$.
  
\end{prop}

\begin{proof}
  We show that $\xi R \xi^{-1} \subset \LO$ if and only if $\xi =
  \diag(\pi^{m_1}, \dots, \pi^{m_n})$.  If $\xi$ is diagonal, the
  result is clear, so we assume that $\xi$ is in Hermite normal form
  and deduce inductively that the off-diagonal entries are zero.

  Let $D = \diag(d_1, \dots, d_n) \in R$, and consider $C = \xi D
  \xi^{-1}$.  We need to examine explicitly the entries of $C$.
  Obviously $\xi D =
  \begin{pmatrix}
    \pi^{m_1}d_1&a_{12}d_2&& \dots && a_{1n}d_n\\
    0&\pi^{m_2}d_2&a_{23}d_3& \dots && a_{2n}d_n\\
    0&0&\ddots&&\\
    0&0&\dots&&\pi^{m_{n-1}}d_{n-1}&a_{n-1\; n}d_n\\
    0&0&&\dots&&\pi^{m_n}d_n
  \end{pmatrix}
  $, and $C_{ij} = \sum_{k=1}^n (\xi D)_{ik} (\xi^{-1})_{kj} =
  \sum_{k=i}^j (\xi D)_{ik} (\xi^{-1})_{kj}$, since both $\xi D$ and
  $\xi^{-1}$ are upper triangular.  Here as is standard
  $(\xi^{-1})_{kj} = (\det \xi)^{-1} (-1)^{k+j} \det \xi(j \slash k)$
  where $\xi(j\slash k)$ is the $(n-1)\times(n-1)$ minor obtained by
  deleting the $j$th row and $k$th column of $\xi$.

  For $1 \le i < n$ we consider the entry $C_{i\; i+1} =
  \sum_{k=i}^{i+1} (\xi D)_{ik} (\xi^{-1})_{k\; i+1}$.  We compute
  \[(\xi^{-1})_{i\; i+1} = (\det \xi)^{-1} (-1)^{2i+1} \det
  \begin{pmatrix}
    \pi^{m_1}&&& \dots &*& \\
    &\ddots&&&&*\\
    &&\pi^{m_{i-1}}&a_{i-1\; i+1}& \dots &&\\
    &&&a_{i\;i+1}&a_{i\; i+2}&\dots \\
    &0&&&\pi^{m_{i+2}}&\\
    &&&0&&\ddots
  \end{pmatrix}= \frac{-a_{i\;i+1}}{\pi^{m_i+m_{i+1}}},
  \]
  so $\ds C_{i\;i+1} =
  (\pi^{m_i}d_i)(\frac{-a_{i\;i+1}}{\pi^{m_i+m_{i+1}}}) +
  (a_{i\;i+1}d_{i+1})(\pi^{-m_{i+1\;i+1}}) =
  \frac{a_{i\;i+1}}{\pi^{m_{i+1}}} (d_{i+1} - d_i).$ Since $C$ must be
  an element of $\LO = M_n(\O)$ we must have $ C_{i\;i+1} \in
  \O$. Since the $d_i$'s are arbitrary we may assume that $\pi \nmid
  (d_{i+1} - d_i)$, so $\ds C_{i\;i+1}=
  \frac{a_{i\;i+1}}{\pi^{m_{i+1}}} (d_{i+1} - d_i) \in \O$ forces
  $a_{i\;i+1} \equiv 0\pmod {\pi^{m_{i+1}}}$.  But we have chosen
  $\xi$ in Hermite normal form which forces $a_{i\;i+1} = 0$.

  Inductively, suppose $a_{ij} = 0$ for $i+1 \le j \le i+\ell$.  We
  show $a_{i\; i+\ell+1} = 0$.  Consider the entry
  \[
  C_{i, i+\ell+1} = \sum_{k=i}^{i+\ell+1} (\xi D)_{ik} (\xi^{-1})_{k\;
    i+\ell+1}= (\xi D)_{i i} (\xi^{-1})_{i,i+\ell+1} (\xi
  D)_{i,i+\ell+1} (\xi^{-1})_{i+\ell+1, i+\ell+1},
  \]
  since $(\xi D)_{i,i+r} = a_{ir}d_r = 0$ for $1\le r\le \ell$. As
  before, there is only one term at issue, $(\xi^{-1})_{i,i+\ell+1} =
  (\det \xi)^{-1} (-1)^{2i+\ell+1} \det \xi(i+\ell+1\slash i)$.  Now
  the minor has the form:
  \[
  \xi(i+\ell+1\slash i)=
  \begin{pmatrix}
    \ddots\\
    &\pi^{m_{i-1}}&a_{i-1,i+1}& \dots\\
    &&a_{i,i+1}&a_{i,i+2}&a_{i,i+3}&\dots&a_{i,i+\ell+1}&\dots\\
    &&\pi^{m_{i+1}}&a_{i+1,i+2}&a_{i+1,i+3}&\dots&a_{i+1,i+\ell+1}&\dots\\
    &&&\pi^{m_{i+2}}&a_{i+2, i+3}&\\
    &&&&\ddots&\ddots\\
    &&&&&\pi^{m_{i+\ell}}&a_{i+\ell,i+\ell+1}\\
    &&&&&&0&\pi^{m_{i+\ell+2}}\\
    &&&&&&&&\ddots\\
  \end{pmatrix}.
  \]
  Recall that by induction, $a_{ij} = 0$ for $i+1 \le j \le i+\ell$.
  As a result, interchanging rows $i, i+1$, then $i+1, i+2$, \dots,
  $i+\ell-1, i+\ell$ produces an upper triangular matrix with
  determinant $\ds \frac{\det \xi}{\pi^{m_i + m_{i+\ell+1}}}
  a_{i,i+\ell+1}$ which because of the interchange of rows differs
  from the determinant of the minor by $(-1)^\ell$.  It now follows
  that
  \[C_{i, i+\ell+1} = (\pi^{m_i}d_i)(-1)\frac{a_{i,i+\ell+1}}{\pi^{m_i
      + m_{i+\ell+1}}} +
  \frac{a_{i,i+\ell+1}d_{i+\ell+1}}{\pi^{m_{i+\ell+1}}} = \frac{a_{i,
      i+\ell+1}}{\pi^{m_{i+\ell+1}}} (d_{i+\ell+1} - d_i).
  \]
  As in the base case, since the $d_k$'s are arbitrary elements of
  $\O$, $\xi$ is in Hermite normal form, and we require $C_{i,
    i+\ell+1} \in \O$, it follows that $a_{i, i+\ell+1} = 0$, which
  completes the proof.
\end{proof}

\begin{cor}\label{cor:characterize-maximal}
  Every maximal order in $M_n(k)$ containing a subring of the form
  $R= \diag(\O, \dots, \O)$ has the form $\Lambda(m_1, \dots, m_n) =
  \begin{pmatrix}
    \O&\p^{m_1 - m_2}&\p^{m_1-m_3}& \dots&\p^{m_1-m_n}\\
    \p^{m_2-m_1}&\O&\p^{m_2-m_3}&\dots&\p^{m_2 - m_n}\\
    \p^{m_3-m_1}&\p^{m_3-m_2}&\ddots&\dots&\p^{m_3 - m_n}\\
    \vdots&\vdots&&\O&\vdots\\
    \p^{m_n - m_1}&\dots&&\p^{m_{n} - m_{n-1}}&\O
  \end{pmatrix}
  $.

  \noindent In particular,  $\Lambda(m_1, \dots, m_n) = \Lambda(0,m_2-m_1,
  \dots, m_n-m_1)$ is the order characterized by $E^{(i,i)}
  \Lambda(m_1, \dots, m_n) E^{(j,j)} = \p^{m_i - m_j} E^{(i,j)}$.
  
\end{cor}

\begin{proof}
  In Proposition~\ref{prop:HNF-to-diagonal}, we observed that the
  maximal orders containing $R$ all have the form $\xi^{-1} M_n(\O)
  \xi$ where $\xi$ is diagonal.  For later convenience in identifying
  vertices with homothety classes of lattices below, we assume that
  $\xi$ has the form $\xi = \diag(\pi^{-m_1}, \dots, \pi^{-m_n})$.
  Thus $\xi^{-1} M_n(\O) \xi$ is certainly contained in the set
  $\Lambda(m_1, \dots, m_n)$.  On the other hand, from
  Proposition~\ref{prop:entries-form-ideals}, it is easily seen that
  the $ij$-entry of $\xi^{-1} M_n(\O) \xi $ is an ideal containing
  $\pi^{m_i-m_j}$, which completes the proof.
\end{proof}

\subsection{Connections to the affine building for
  $SL_n(k)$}\label{sec:geometric-connections}

To introduce the connection between split orders in $B$ and convex
polytopes in affine buildings requires a bit of background which we
present here in abbreviated form; the books by Brown \cite{Brown} and
Garrett \cite{Garrett} are two excellent resources for further
details.  Classically, affine buildings are associated to $p$-adic
groups, e.g., $SL_n(k)$, and are characterized as simplicial complexes
whose simplicial structure is determined by subgroups and cosets of
the $p$-adic group being studied.  Here, we give a well-known but more
arithmetic characterization.  To present the standard nomenclature,
the simplicial complex which is the building is itself the union of
subcomplexes called apartments, all of which are isomorphic.
Apartments of an affine building are tilings of Euclidean space,
and the structure of the tiling is determined by the associated
Coxeter diagram which encodes the generators and relations of the Weyl
group associated to the $p$-adic group.

The affine building for $SL_n(k)$ is an $(n-1)$-dimensional simplicial
complex in which the maximal orders in $B=M_n(k)$ comprise the
vertices.  Apartments in the building are $(n-1)$-complexes, whose
structure is captured by a tessellation of $\R^{n-1}$. We give a
concrete realization; see \cite{Brown} or \cite{Garrett} for further
details.  Let $V$ be an $n$-dimensional vector space over the local
field $k$, and identify $B = M_n(k)$ with $End_k(V)$.  Let $L$ be any
lattice (free $\O$-module of rank $n$) in $V$.  The
homothety class of $L$, denoted $[L]$, is simply the set
of lattices $\{\lambda L\mid \lambda \in k^\times\}$.

It is easy to check that for two lattices $L$ and $M$, the homothety
classes $[L] = [M]$ iff $End_\O(L) = End_\O(M)$, and that as $L$ runs
through the set of lattices of $V$, $End_\O(L)$ runs through the set
of maximal orders of $B$.  Thus, the vertices of our
building originally given by maximal orders in $B$, may instead be
identified with the homothety classes of lattices in $V$.  To introduce
the simplicial structure, we define the notion of incidence: we say
that two vertices are incident if there are lattices $L$ and $L'$
representing the vertices such that $\pi L \subseteq L' \subseteq L$.
Note in this case, $\pi L' \subseteq \pi L \subseteq L'$, so the
definition of incidence is symmetric, and defines the edges
(1-simplicies) in the building.  An $m$-simplex is characterized by
lattices $L_i$ (representing its vertices) satisfying $\pi L_0
\subsetneq L_1 \subsetneq \cdots \subsetneq L_m \subsetneq L_0$, or
equivalently flags of length $m$ in the $\O/\pi \O$-vector space
$L_0/\pi L_0$.  The maximal simplicies ($(n-1)$-simplicies) are called
the chambers of the building.

To make things even more concrete, we note (\cite{Garrett}) that there
is a one-to-one correspondence between sets of $n$ linearly
independent lines in $V$ (frames) and apartments in the building for
$SL_n(k)$. In particular, every vertex in a fixed apartment can be
represented by a lattice of the form $\O \pi^{\nu_1}e_1 \oplus \cdots
\oplus \O \pi^{\nu_n} e_n$ for some fixed basis $\{e_1, \dots, e_n\}$
of $V$, and where the $\nu_i$ range over all elements of $\Z$.  Since each
vertex in the apartment is the homothety class of a lattice $\O
\pi^{\nu_1}e_1 \oplus \cdots \oplus \O \pi^{\nu_n} e_n$, we may simply
identify the vertices in an apartment in the $SL_n(k)$ building with
the elements of $\Z^n/\Z(1,1,\dots,1)$, where we represent the
homothety class of $\O \pi^{\nu_1}e_1 \oplus \cdots \oplus \O
\pi^{\nu_n} e_n$ by $[\nu_1, \dots, \nu_n]$ or after normalizing, by
$[0, \nu_2 - \nu_1, \dots, \nu_n-\nu_1]$.

To recast some of our earlier algebraic results in this geometric
setting, we let $V$ be as above, fix a basis $\{e_1, \dots, e_n\}$ for
$V$ and let $L_0$ be the $\O$-lattice with basis $\{e_i\}$.
Identifying $\End_\O(L_0)$ with $\LO = M_n(\O)$, we observe that for
$\xi \in B^\times$, $\xi^{-1} \LO \xi = \End(\xi^{-1} L_0)$, so all
maximal orders in $B$ have the form $\End(\xi^{-1} L_0)$ for some $\xi
\in B^\times$.  In Corollary~\ref{cor:characterize-maximal}, we showed
that every maximal order containing $R = \diag(\O,\dots,\O)$ can be
expressed as $\Lambda(m_1, \dots, m_n) = \Lambda(0,m_2-m_1, \dots,
m_n-m_1)$.  So taking $\xi=\diag(1, \pi^{-m_2}, \dots, \pi^{-m_n})$,
we can identify $\Lambda(0,m_2, \dots, m_n)$ with the homothety class
of the lattice $\xi^{-1}L_0 = \O e_1 \oplus \O\pi^{m_2} e_2 \oplus
\dots \oplus \O \pi^{m_n}e_n$ which we denote $[0, m_2, \dots, m_n]$.
Thus the set of maximal orders containing $R$ can be represented as
vertices of the building given by homothety classes $[0, m_2, \dots,
m_n]$, $m_i \in \Z$.

\begin{remark}\label{remark:restrict-to-apt}
  The significance of the above characterization is twofold.  First,
  every maximal order in this fixed apartment contains $R$, so
  that all such maximal orders are split orders.  More significantly
  is that if we wish to consider orders $S$ which contain $R$, the set
  of maximal orders which contain $S$ all lie in a given apartment.
  Of course there may be many such apartments, but the ability to
  restrict to a fixed apartment leads not only to the concrete
  algebraic representation, but more importantly to the geometric one
  we develop below.
\end{remark}

\section{Geometric considerations}

Our goal is to give a geometric characterization of split orders, and
we begin in our restricted setting of split orders $S$ of $B = M_n(k)$
with $R = \diag(\O, \dots, \O)\subset S \subset B$.  By
Remark~\ref{remark:restrict-to-apt}, we can and do fix an apartment
$\A_0$ which contains all the maximal orders $\Lambda(0,m_2,\dots,m_n)$
that contain a given $S$.  Via a fixed basis for $V$ (which yields the
frame defining $\A_0$), we identify the apartment with $\R^{n-1}
\cong\{0\} \times \R^{n-1} \subset \R^{n}$; the set of vertices in
$\A_0$ is identified with $\{0\} \times \Z^{n-1} \cong \Z^n/\Z(1, \dots,
1)$.  As noted after Proposition~\ref{prop:entries-form-ideals}, we
adopt the succinct presentation of $S$ as $S = (\p^{\nu_{ij}}) =
\left(\begin{smallmatrix}
    \O&&\p^{\nu_{ij}}\\
    &\ddots&\\
    \p^{\nu_{ij}}&&\O\\
        \end{smallmatrix}\right)$, with $\nu_{ij} \in \Z$, $\nu_{ii} =
      0$.

For a maximal order
$\Lambda(0, m_2, \dots, m_n) =
\begin{pmatrix}
\O&\p^{-m_2}&\p^{-m_3}& \dots&\p^{-m_n}\\
\p^{m_2}&\O&\p^{m_2-m_3}&\dots&\p^{m_2 - m_n}\\  
\p^{m_3}&\p^{m_3-m_2}&\ddots&\dots&\p^{m_3 - m_n}\\
\vdots&\vdots&&\O&\vdots\\
\p^{m_n}&\dots&&\p^{m_{n} - m_{n-1}}&\O
\end{pmatrix}
$, we have 
$S \subset \Lambda(0, m_2, \dots, m_n)$ if and only if (setting $m_1 = 0$)
\begin{equation}\label{eqn:convex-hull}
 -\nu_{ji} \le m_i - m_j \le \nu_{ij} \mbox{ for all } i, j.
\end{equation}

Given our identification of the apartment $\A_0$ with $\R^{n-1}$, the
equations of the form $L_{ij} := x_i - x_j = \nu \in \Z$ are hyperplanes in
$\R^{n-1}$ and represent a subset of the walls in the apartment; they
represent all of the walls if $n=2,3$.  It is clear that the
inequalities 
\begin{equation}
\label{eqn:linear-inequalities}
 -\nu_{ji} \le L_{ij} = x_i - x_j \le \nu_{ij}
\end{equation}
define a convex polytope
in $\R^{n-1}$ which we denote by $C_S$.

The immediate aim of this section is to establish a one-to-one
correspondence between split orders containing $R$ and convex
polytopes of this form in the apartment $\A_0$.  We have already seen
(\ref{eqn:convex-hull}) that the convex hull determined by the walls
of the building which contain the set of maximal orders containing a
given split order forms a convex polytope.  We now further show that
the split order is the intersection of the maximal orders contained in
that polytope.

\begin{definition}
  Given our fixed apartment $\A_0$, let $\C$ denote the set of convex
  polytopes determined by systems of inequalities as in
  (\ref{eqn:linear-inequalities}); we denote a typical element in $\C$
  as $C(\bnu)$, $\bnu = (\nu_{ij}) \in M_n(\Z)$.  We shall require
  that $\bnu$ (or $C(\bnu)$) be \textit{reduced}, meaning the convex
  region determined by the inequalities
  (\ref{eqn:linear-inequalities}) contain at least one vertex of the
  building, and each of the hyperplanes determined by the $\nu_{ij}$
  meets the convex region. In the usual terminology of convex
  geometry, each of the given hyperplanes $L_{ij} = \nu_{ij}$ or
  $L_{ij} = -\nu_{ji}$ is a supporting hyperplane.
\end{definition}

\begin{remark}
  Note that since $x_1 = 0$ in our characterization of the
  apartment $\A_0$, the inequalities $-\nu_{1i} \le x_i - x_1 \le \nu_{i1}$
  reduce to  $-\nu_{1i} \le x_i \le \nu_{i1}$, so that $C(\bnu)$
  always defines a compact convex region, hence one containing only
  finitely many vertices.
\end{remark}

\begin{prop}\label{prop:first-half}
  Let $C=C(\bnu) \in \C$, and let $S_C = (\p^{\mu_{ij}}) =
\bigcap_{\Lambda \in C} \Lambda$ be the split order which is the
intersection of all maximal orders in $C(\bnu)$. Then
$\mu_{ij} = \nu_{ij}$ for all $i,j$.
\end{prop}

\begin{proof}
  Let $\Lambda_k$ index the maximal orders (vertices) in $C(\bnu)$,
  and denote $\Lambda_k = (\p^{\lambda_{ij}^{(k)}})$.  Since $S_C=
  (\p^{\mu_{ij}})$ is the intersection of the $\Lambda_k$, it is clear
  that $\mu_{ij} = \max_k\{ \lambda_{ij}^{(k)}\}$, so $\mu_{ii} =
  \lambda_{ii}^{(k)} = 0$.  For each $i < j$ we have $-\nu_{ji} \le
  \lambda_{ij}^{(k)}\le \nu_{ij}$, so $\mu_{ij} = \max_k\{
  \lambda_{ij}^{(k)}\}$ and $\mu_{ji} = \max_k\{
  \lambda_{ji}^{(k)}\}=\max_k\{ -\lambda_{ij}^{(k)}\} =
  -\min_k\{\lambda_{ij}^{(k)}\}$.  However, since for the convex
  region $C(\bnu)$, we require that $\bnu$ be reduced, there are
  maximal orders on the boundary of the region achieving each of the
  bounding limits.  Thus for $i < j$, $\mu_{ij} = \max_k\{
  \lambda_{ij}^{(k)}\} = \nu_{ij}$, while $\mu_{ji} =
  -\min_k\{\lambda_{ij}^{(k)}\} = \nu_{ji}$.
\end{proof}

Let's examine the correspondence as it now stands.  Given $C = C(\bnu)
\in \C$, we form $S_C = (\p^{\mu_{ij}}) = \bigcap_{\Lambda \in C}
\Lambda$, and since $\mu_{ij} = \nu_{ij}$, we have $C(\bmu)
= C(\bnu)$ which is half of the desired correspondence between split
orders and convex polytopes. Perhaps more succinctly we have:
\[
C = C(\bnu) \mapsto S_C = (\p^{\mu_{ij}}) = \bigcap_{\Lambda \in C}
\Lambda \mapsto C(\bmu) = C.
\]

To establish the other half of the correspondence, 
\[S=(\p^{\nu_{ij}})\mapsto C(\bnu) \mapsto \bigcap_{\Lambda \in C(\bnu)}
\Lambda = (\p^{\mu_{ij}})= S,
\]
significantly more
effort is required.  Consider a subset of $M_n(k)$ having the form $S
= (\p^{\nu_{ij}})$.  A necessary condition that $S$ be contained in
some maximal order is that $\nu_{ij} + \nu_{ji} \ge 0$ for all $i,
j$. Given that necessary condition, $S$ determines a convex polytope
$C_S = C(\bnu)$ via the pairs of inequalities in
(\ref{eqn:linear-inequalities}).  The potential difficulty is that
different subsets $S$ can determine the same convex region.  The
following example demonstrates the difficulty and suggests its
resolution.

\begin{example}
  Consider $S = \pmatt{\O}{\O}{\p}{\p^3}{\O}{\p}{\p^3}{\p^2}{\O}$ and
  $S' = \pmatt{\O}{\O}{\p^2}{\p^3}{\O}{\p}{\p^3}{\p^2}{\O}$.  The
  diagram below (points have coordinates $[0,x_2,x_3]$) illustrates
  that $S$ and $S'$ determine the same convex region via the
  inequalities (\ref{eqn:linear-inequalities}):
    \begin{align*}
    0 &\le x_2 \le 3 \hfill &0 \le x_2 \le 3\\
   S:\quad -1&\le x_3 \le 3 \hfill &S':\qquad-2\le x_3 \le 3\\
    -1 &\le x_3 - x_2 \le 2 \hfill &-1 \le x_3 - x_2 \le 2\\
    \end{align*}

    From the diagram, we see that the hyperplane $x_3 = -2$ does not
    intersect the convex polytope, while the hyperplane $x_3 = -1$
    does so in precisely one point, though neither is actually required to
    determine the convex region.
    
\includegraphics[height=3in]{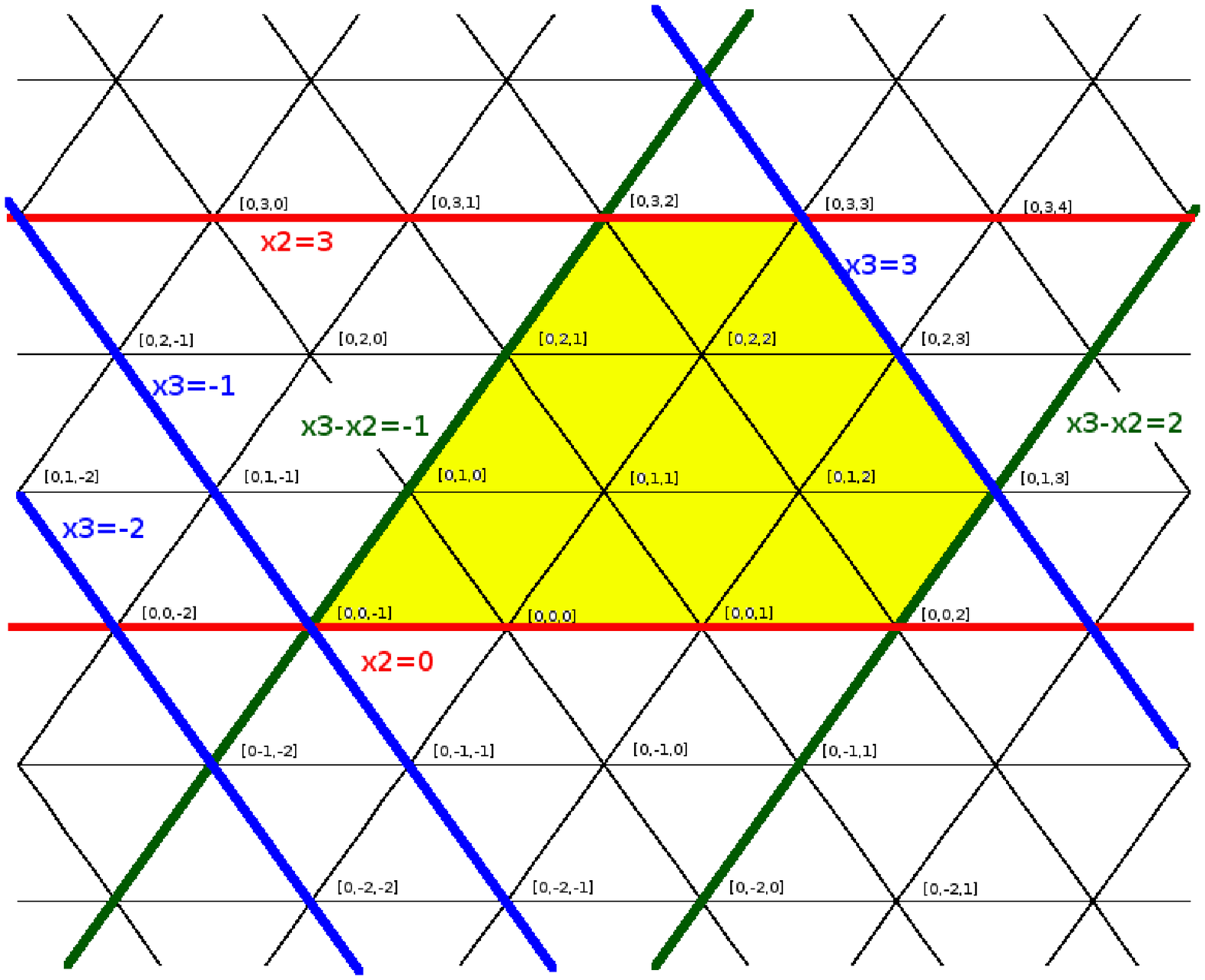}

 $S$ is an order since it is the intersection of the maximal orders in
    this convex region, however $S'$ is not.  Indeed, $S$ is the
    intersection of those maximal orders on the boundary of the convex
    polytope which determine it:
    \begin{align*}
    S &= \Lambda(0,0,-1)\cap\Lambda(0,3,2)\cap\Lambda(0,3,3) \cap
    \Lambda(0,1,3)\cap \Lambda(0,0,2)\\
    &=\Lambda(0,0,-1)\cap\Lambda(0,3,3)
    \cap \Lambda(0,0,2)\\
    &= \Lambda(0,0,-1)\cap\Lambda(0,3,2)
    \cap\Lambda(0,1,3) \\
    &=\spmatt{\O}{\O}{\p}{\O}{\O}{\p}{\p^{-1}}{\p^{-1}}{\O}\cap
    \spmatt{\O}{\p^{-3}}{\p^{-2}}{\p^3}{\O}{\p}{\p^{2}}{\p^{-1}}{\O}\cap
    \spmatt{\O}{\p^{-1}}{\p^{-3}}{\p}{\O}{\p^{-2}}{\p^{3}}{\p^{2}}{\O}
  \end{align*}
  
  On the other hand it is easy to see that $S'$ is not an order.  By
  Proposition~\ref{prop:entries-form-ideals}, a necessary condition
  that a subset $S = (\p^{\nu_{ij}}) \supset R$ be an order is that it
  be closed under multiplication, which requires $\nu_{ik} + \nu_{kj}
  \ge \nu_{ij}$, $\nu_{kk} = 0$ for all $i, j, k$.  We note that in
  $S'$, $\nu_{12} + \nu_{23} = 0 + 1 \not\ge 2 = \nu_{13}$.
\end{example}

The key to establishing the other half of the desired correspondence is to
connect the failure to be an order with a geometric condition.  This
leads us to the the following definition.

 \begin{definition}
   Let $S = (\p^{\nu_{ij}})$ be a subset of $M_n(k)$ which contains
   $R$ and satisfies $\nu_{ij}+\nu_{ji} \ge 0$ for all $i, j$.  Call
   $S$ \textit{reduced} if the convex region it determines, $C(\bnu)$,
   is reduced.
  \end{definition}

  \begin{prop}\label{prop:second-half} Let $S = (\p^{\nu_{ij}})$ be as above.  Then $S$ is an
    order if and only if $S$ is reduced.
  \end{prop}
 
\begin{remark}
  Note that if $S$ is not reduced, there is an $\overline S \supset S$
  which is reduced (hence an order), and which determines exactly the same
  convex polytope.
  \end{remark}

\begin{proof}
  One direction is quite easy.  If $S$ is reduced, the bounds on the
  inequalities defining the convex polytope
  (\ref{eqn:linear-inequalities}) are sharp, and from the arguments
  above, the intersection of all the maximal orders in that convex
  polytope equals $S$, that is $S$ is the intersection of the maximal
  orders containing it, hence $S$ is an order.

  Note that by Proposition~\ref{prop:entries-form-ideals}, $S =
  (\p^{\nu_{ij}})$ is an order if and only if $\nu_{ik} + \nu_{kj} \ge
  \nu_{ij}$ for every $i, j, k$.  So to establish the converse of our
  theorem, we show that $\nu_{ik} + \nu_{kj} \ge \nu_{ij}$ for every
  $i, j, k$ implies $\bnu = (\nu_{ij})$ (i.e., $C(\bnu)$) is reduced.  We
  proceed by contradiction, so we assume that there exist $i_0, j_0$
  such that $x_{i_0} - x_{j_0} = \nu_{i_0j_0}$ does not intersect
  $C(\bnu)$.
  
Since $x_1 = 0$, there is some asymmetry in the expression $x_{i_0} -
x_{j_0}$ when one of $i_0, j_0 = 1$, so we separate the proof into cases
beginning with the generic case.

\noindent\textbf{ Case: $\bm{i_0, j_0 \ne 1}$.}  If the hyperplane $x_{i_0} -
x_{j_0} = \nu_{i_0j_0}$ does not intersect $C(\bnu)$, we have $x_{i_0}
- x_{j_0} < \nu_{i_0j_0}$ for all $(x_i) \in C(\bnu)$. Note that the
symmetric case $x_{i_0} - x_{j_0} > -\nu_{j_0i_0}$ is equivalent to
$x_{j_0} - x_{i_0} < \nu_{j_0i_0}$ so we consider only $x_{i_0} -
x_{j_0} < \nu_{i_0j_0}$.  Let $\bb = (b_i) \in C(\bnu)$ achieve a
maximum for $x_{i_0} - x_{j_0}$, say $b_{i_0} - b_{j_0} = \mu_{i_0j_0}
< \nu_{i_0j_0}$.  To arrive at the desired contradiction, we use $\bb$
to construct a point $\bb'= (b'_i) \in C(\bnu)$ with $\mu_{i_0j_0} <
b'_{i_0} - b'_{j_0} \le \nu_{i_0j_0}$.  Note, throughout the proof we
use without further mention that all hyperplanes have the form $x_i -
x_j = \nu \in \Z$.

We set a bit of notation.  Let $\be_\ell$ be the $\ell$th standard
basis vector in $\R^n$, and for $k \ne 1, i_0, j_0$, let
\begin{equation*}
  \alpha_k=
  \begin{cases}
    1&\textrm{if }b_k - b_{j_0} = \nu_{kj_0},\\0&\textrm{otherwise},
  \end{cases}\hspace{.5in}
  \beta_k = 
  \begin{cases}
    1&\textrm{if }b_{i_0} - b_k = \nu_{i_0k},\\0&\textrm{otherwise}.
  \end{cases}
\end{equation*}
To define $\bb'$ we need to increase the difference
$b_{i_0} - b_{j_0}$, either by increasing $b_{i_0}$ or decreasing
$b_{j_0}$ and adjust the other coordinates to satisfy all the
remaining convexity bounds. We put
\begin{equation}\label{eqn:bprime-case1}
  \bb' =
  \begin{cases}
    \bb - \be_{j_0} - \sum_{k\ne 1, i_0, j_0} \alpha_k \be_k&
    \textrm{if } b_{i_0} = \nu_{i_01},\\
    \bb + \be_{i_0} +  \sum_{k\ne 1, i_0, j_0} \beta_k \be_k&\textrm{if } b_{i_0} < \nu_{i_01}.
  \end{cases}
\end{equation}

\noindent Subcase A. We begin with the case where $b_{i_0} = \nu_{i_01}$ and 
$\bb' = \bb - \be_{j_0} - \sum_{k\ne 1, i_0, j_0} \alpha_k \be_k$.

First we show that $-\nu_{1i} \le b'_i = b'_i - b'_1 \le \nu_{i1}$ for
all $i$.  This is clear for $b'_1 = 0$ and  $b'_{i_0}= b_{i_0} =
\nu_{i_01}$.  We note $b'_{j_0} = b_{j_0} - 1\le \nu_{j_01} - 1\le
\nu_{j_01}$.  To see $b'_{j_0} \ge -\nu_{1j_0}$, note that 
$b_{i_0} - b_{j_0} = \mu_{i_0j_0} < \nu_{i_0j_0} \le \nu_{i_01} +
\nu_{1j_0}$ by assumptions on $S=(\p^{\nu_{ij}})$ and $\bb$, so 
\begin{equation}\label{eqn:bj0}
  b_{j_0} = b_{i_0} - \mu_{i_0j_0} = \nu_{i_01} - \mu_{i_0j_0} >
  \nu_{i_01} - \nu_{i_0j_0} \ge -\nu_{1j_0}.
\end{equation}
Thus $b_{j_0} >  -\nu_{1j_0}$ implies $b'_{j_0} = b_{j_0} - 1 \ge
-\nu_{1j_0}$ as desired.  Next we finish the remaining inequalities of
the form $-\nu_{1k} \le b'_k = b'_k -b'_1 \le \nu_{k1}$, $k \ne 1, i_0, j_0$.

By the definition of $\bb'$, if $b_k - b_{j_0} < \nu_{kj_0}$, then
$b_k' = b_k$, so there is no issue.  If $b_k - b_{j_0} = \nu_{kj_0}$,
then $b_k' = b_k - 1$, so of course $b'_k \le \nu_{k1}$.  To see $b'_k
\ge -\nu_{1k}$, we suppose not, so $b'_k = b_k - 1 < -\nu_{1k}$, hence
$b_k < -\nu_{1k} + 1$.  On the other hand, $\bb \in C(\bnu)$ implies
$b_k \ge -\nu_{1k}$ from which we deduce $b_k = -\nu_{1k}$.  Now $b_k
- b_{j_0} = \nu_{kj_0}$ implies $b_{j_0} = b_k - \nu_{kj_0} =
-\nu_{1k} - \nu_{kj_0} \le -\nu_{1j_0}$ contrary to equation~(\ref{eqn:bj0}).

Next we must consider bounds on $b'_k - b'_\ell$ where $k, \ell \ne 1$,
and where $\{k,\ell\} \cap \{i_0,j_0\}$ has cardinality 0, 1, or 2.

First observe that since $\mu_{i_0j_0} < \nu_{i_0j_0}$ (and both are
  integers), 
\begin{equation*}
  -\nu_{j_0i_0} \le b_{i_0} - b_{j_0} <  b'_{i_0} - b'_{j_0} = b_{i_0}
  - (b_{j_0} - 1) = \mu_{i_0j_0} + 1 \le \nu_{i_0j_0}
\end{equation*}

Next consider
\begin{equation*}
-\nu_{j_0k} \le b_k - b_{j_0} \le  b'_k - b'_{j_0} =
\begin{cases}
  b_k - b_{j_0} + 1&\textrm{if }b_k - b_{j_0} < \nu_{kj_0},\\
   b_k - b_{j_0}&\textrm{if }b_k - b_{j_0} = \nu_{kj_0}
\end{cases}
\qquad\le \nu_{kj_0}.
\end{equation*}
Similarly, since
\begin{equation*}
  b'_{i_0} - b'_k =
  \begin{cases}
    b_{i_0} - b_k + 1&\textrm{if }b_k - b_{j_0} = \nu_{kj_0},\\
    b_{i_0} - b_k &\textrm{if }b_k - b_{j_0} < \nu_{kj_0},\\
  \end{cases}
\end{equation*}
it is clear that $-\nu_{ki_0} \le b'_{i_0} - b'_k$, and the only
potential issue with the upper bound is when $b_k - b_{j_0} =
\nu_{kj_0}$.  If indeed $ b'_{i_0} - b'_k > \nu_{i_0k}$, then $b_{i_0}
- b_k > \nu_{i_0k} - 1$ which means $b_{i_0} - b_k = \nu_{i_0k}$.  But
this together with $b_k - b_{j_0} = \nu_{kj_0}$ implies $b_{i_0} -
b_{j_0} = \nu_{i_0k} + \nu_{kj_0} \ge \nu_{i_0j_0}$, but by hypothesis
$b_{i_0} - b_{j_0} = \mu_{i_0j_0} < \nu_{i_0j_0}$, a contradiction.

Finally, we come to the case $b'_k - b'_\ell$ where $k, \ell \ne 1$,
and where $\{k,\ell\} \cap \{i_0,j_0\} = \emptyset$.  We see that
\begin{equation*}
  b'_k - b'_\ell =
  \begin{cases}
    b_k - b_\ell + 1&\textrm{if }b_k - b_{j_0} < \nu_{kj_0} \textrm{
      and } b_\ell - b_{j_0} = \nu_{\ell j_0},\\
     b_k - b_\ell - 1&\textrm{if }b_k - b_{j_0} = \nu_{kj_0} \textrm{
      and } b_\ell - b_{j_0} < \nu_{\ell j_0},\\
     b_k - b_\ell&\textrm{otherwise.}
  \end{cases}
\end{equation*}
Everything is clear except for the upper bound in the first case and
the lower bound in the second case. Assuming $b_k - b_{j_0} < \nu_{kj_0}$
      and $b_\ell - b_{j_0} = \nu_{\ell j_0}$, if $b_k - b_\ell + 1 >
      \nu_{k\ell}$, then  $b_k - b_\ell = \nu_{k\ell}$.  This implies
      $b_k - b_{j_0} = \nu_{k\ell} + \nu_{\ell j_0} \ge \nu_{k j_0}$,
      a contradiction.  Analogously, assuming $b_k - b_{j_0} =
      \nu_{kj_0}$  and $b_\ell - b_{j_0} < \nu_{\ell j_0}$, if $ b_k -
      b_\ell - 1 < -\nu_{\ell k}$ then $b_k - b_\ell = -\nu_{\ell k}$
      which $b_\ell - b_{j_0} = \nu_{\ell k} + \nu_{k j_0} \ge
      \nu_{\ell j_0}$, a contradiction.

      \noindent Subcase B.  Here we assume $b_{i_0} < \nu_{i_01}$ and
      $\bb' = \bb + \be_{i_0} + \sum_{k\ne 1, i_0, j_0} \beta_k
      \be_k$, \newline $\beta_k =
  \begin{cases}
    1&\textrm{if }b_{i_0} - b_k = \nu_{i_0k},\\0&\textrm{otherwise}.
  \end{cases}$

  First we show that $-\nu_{1i} \le b'_i = b'_i - b'_1 \le \nu_{i1}$
  for all $i$.  This is clear for $b'_1 = 0$ and $b'_{j_0}= b_{j_0}$.
  We note $-\nu_{1i_0} \le b_{i_0} < b'_{i_0} = b_{i_0} + 1 \le
  \nu_{i_0 1}$ since  $b_{i_0} < \nu_{i_01}$.  For $k\ne 1, i_0, j_0$,
  the only issue is when $b_{i_0} - b_k = \nu_{i_0 k}$ in which case
  $b'_k = b_k + 1$, and then only concerns the upper bound.  If
  $b'_k > \nu_{k1}$, then $b_k = \nu_{k1}$, so that $b_{i_0} = b_k +
  \nu_{i_0 k} =  \nu_{i_0 k} + \nu_{k1} \ge \nu_{i_0 1}$, contrary to
  assumption.

  Next observe
  \begin{equation*}
    -\nu_{j_0i_0} \le b_{i_0} - b_{j_0} < b'_{i_0} - b'_{j_0} = b_{i_0}
    - b_{j_0} + 1 = \mu_{i_0 j_0} + 1 \le \nu_{i_0j_0}.
  \end{equation*}

We also have
\begin{equation*}
  -\nu_{ki_0} \le b'_{i_0} - b'_k =
  \begin{cases}
    b_{i_0} - b_k&\textrm{if } b_{i_0} - b_k = \nu_{i_0k},\\
     b_{i_0} - b_k+ 1&\textrm{if } b_{i_0} - b_k < \nu_{i_0k}\\
  \end{cases}
\qquad \le \nu_{i_0 k}.
\end{equation*}

Similarly, since
\begin{equation*}
  b'_k - b'_{j_0} =
  \begin{cases}
    b_k - b_{j_0}&\textrm{if } b_{i_0} - b_k < \nu_{i_0k},\\
     b_k - b_{j_0} + 1&\textrm{if } b_{i_0} - b_k = \nu_{i_0k},\\
  \end{cases}
\end{equation*}
the only issue is with the upper bound when $ b_{i_0} - b_k =
\nu_{i_0k}$.  If $b'_k - b'_{j_0}> \nu_{kj_0}$, then $b_k - b_{j_0} =
\nu_{kj_0}$.  This together with $ b_{i_0} - b_k =
\nu_{i_0k}$ implies $b_{i_0} - b_{j_0} = \nu_{i_0k} + \nu_{kj_0} \ge
\nu_{i_0j_0}$, contrary to our original assumption.

Finally, we come to the case $b'_k - b'_\ell$ where $k, \ell \ne 1$,
and where $\{k,\ell\} \cap \{i_0,j_0\} = \emptyset$.  We see that
\begin{equation*}
  b'_k - b'_\ell =
  \begin{cases}
    b_k - b_\ell + 1&\textrm{if }b_{i_0} - b_{k} = \nu_{i_0k} \textrm{
      and } b_{i_0} - b_{\ell} < \nu_{i_0\ell},\\
     b_k - b_\ell - 1&\textrm{if }b_{i_0} - b_{k} < \nu_{i_0k} \textrm{
      and } b_{i_0} - b_{\ell} = \nu_{i_0\ell},\\
     b_k - b_\ell&\textrm{otherwise.}
  \end{cases}
\end{equation*}
Everything is clear except for the upper bound in the first case and
the lower bound in the second case.  Assuming $b_{i_0} - b_{k} =
\nu_{i_0k}$ and $b_{i_0} - b_{\ell} < \nu_{i_0\ell},$ if $ b_k -
b_\ell + 1 > \nu_{k\ell}$, then $b_k - b_\ell = \nu_{k\ell}$, so that
$b_{i_0} - b_\ell = \nu_{i_0k} + \nu_{k\ell}\ge \nu_{i_0\ell}$,
contrary to assumption.  Analogously, assuming $b_{i_0} - b_{k} <
\nu_{i_0k}$ and $b_{i_0} - b_{\ell} = \nu_{i_0\ell},$ if $b_k - b_\ell
- 1< \nu_{\ell k}$, then $b_k - b_\ell = -\nu_{\ell k}$, so $b_{i_0} -
b_k = \nu_{i_0\ell} + \nu_{\ell k} \ge \nu_{i_0 k}$, contrary to assumption.

\noindent\textbf{ Case: $\bm{i_0 \textrm{ or } j_0 = 1}$.} We choose
$\bb$ as in the first case and define $\alpha_k$ and $\beta_k$ exactly
as before (noting the obvious redundant conditions on $k$).  We put
\begin{equation}\label{eqn:bprime-case2}
  \bb' =
  \begin{cases}
    \bb - \be_{j_0} - \sum_{k\ne 1, i_0, j_0} \alpha_k \be_k&
    \textrm{if } i_0 = 1,\\
    \bb + \be_{i_0} +  \sum_{k\ne 1, i_0, j_0} \beta_k \be_k&\textrm{if } j_0 = 1.
  \end{cases}
\end{equation}
Then these boundary cases are handled in exactly the same way as above
with no further insights required, and this completes the proof of the
theorem.
\end{proof}

We summarize both pieces of the correspondence as

\begin{thm}
  There is a one-to-one correspondence between convex polytopes in
  $\C$ determined by the walls of the apartment $\A_0$ in the building
  for $SL_n(k)$ and split orders in $M_n(k)$ which contain $R$.  The
  maps $C = C(\bnu) \mapsto S_C = (\p^{\mu_{ij}}) = \bigcap_{\Lambda
    \in C} \Lambda $ and $S=(\p^{\nu_{ij}})\mapsto C(\bnu)$ are
  inverse to one another.
\end{thm}

\begin{proof}
  Given $C(\bnu) \in \C$, we have $\bnu$ is reduced, so by
  Proposition~\ref{prop:first-half}, $C(\bnu) \mapsto S_C =
  \bigcap_{\Lambda \in C} \Lambda = (\p^{\nu_{ij}}) \mapsto C(\bnu)$.
  On the other hand if $S = (\p^{\mu_{ij}})$ is a split order, then by
  Proposition~\ref{prop:second-half}, $\bmu = (\mu_{ij})$ is reduced,
  so that $S = (\p^{\mu_{ij}}) \mapsto C(\bmu) \mapsto
  \bigcap_{\Lambda \in C(\bmu)} \Lambda = (\p^{\mu_{ij}})$ by
  Proposition~\ref{prop:first-half}.
\end{proof}

Now we generalize the above results to that of our general notion of a split
order.  We begin by showing the intersection of any finite collection
of maximal orders in a fixed apartment is a split order.

\begin{prop}\label{prop:intersection-is-split}
  Let $\A$ be any apartment in the affine building for $SL_n(k)$, and
  let $\Lambda_1,\dots, \Lambda_r$ be maximal orders in $M_n(k)$
  corresponding to vertices in $\A$.  Then $S = \bigcap_{i=1}^r
  \Lambda_i$ is a split order.
\end{prop}

\begin{proof}
  Our original fixed apartment $\A_0$ corresponds to the basis
  $\{e_i\}$ of the vector space $V$.  Let $\{f_i\}$ be a basis of $V$
  whose frame determines the apartment $\A$.  Let $\gamma \in GL_n(k)$
  be the change of basis matrix taking $e_i$ to $f_i$.  Each maximal
  order $\Lambda_k = \End_\O(L_k)$ for a lattice $L_k = \oplus \O
  \pi^{a_i^{(k)}} f_i$. Let $\tilde L_k = \gamma^{-1} L_k = \oplus \O
  \pi^{a_i^{(k)}} e_i$ and $\tilde \Lambda_k = \End_\O(\tilde L_k)$.
  Then
\[\Lambda_k = \End_\O(L_k) = \End(\gamma\tilde L_k) =
\gamma \End_\O(\tilde L_k) \gamma^{-1} = \gamma \tilde\Lambda_k \gamma^{-1}.
\]
Now all of the $\tilde\Lambda_k$ are maximal orders in $\A_0$, which by
Remark~\ref{remark:restrict-to-apt} all contain $R$.  Thus, 
$S = \bigcap_{i=1}^r \Lambda_i \supset \gamma R \gamma^{-1}$, hence is a
split order. 
\end{proof}

Next we consider the converse.

\begin{prop}\label{prop:split-is-intersection}
  Suppose that $S$ is an order of $B = M_n(k)$ which contains $\gamma
  R \gamma^{-1}$ for some $\gamma \in B^\times$.  Then $S$ is the
  intersection of maximal orders lying in a convex polytope in the
  apartment $\A = \gamma \A_0$.
\end{prop}

\begin{proof}
  If $S \supset \gamma R \gamma^{-1}$, then $\gamma^{-1} S \gamma$ is
  an order of $B$ containing $R$.  By
  Propositions~\ref{prop:first-half} and \ref{prop:second-half}, 
$\gamma^{-1}S\gamma = (\p^{\bnu}) = \bigcap_{\tilde\Lambda \in C(\bnu)}
\tilde\Lambda$, that is $\bnu$ is reduced and $\gamma^{-1}S\gamma$ is
the intersection of all the maximal orders $\tilde\Lambda$ in the convex
polytope $C(\bnu)$.  It follows that 
\[S = \gamma \left(\bigcap_{\tilde\Lambda \in C(\bnu)}
\tilde\Lambda\right) \gamma^{-1} =  \bigcap_{\tilde\Lambda \in C(\bnu)} \gamma
\tilde \Lambda \gamma^{-1}.
\]
Now let $\tilde\Lambda = \End_\O(\tilde L)$ and $\tilde\Lambda'
= \End_\O(\tilde L')$ be two maximal orders in $\C(\bnu)$.  Then
$\gamma \tilde \Lambda \gamma^{-1} = \End_\O(\gamma \tilde L)$ and
$\gamma \tilde \Lambda' \gamma^{-1} = \End_\O(\gamma \tilde L')$.
Since $\gamma$ can simply be viewed as a change of basis matrix, the
elementary divisors of $L'$ in $L$, denoted $\{L~:~L'\}$, equal those
of $\gamma L'$ in $\gamma L$, that is $\{L:L'\} = \{\gamma
L:\gamma L'\}$.  Moreover, since the incidence relations among vertices in
the building are determined by chains of lattices whose relative
containments in an apartment are completely determined by the
elementary divisors, we see that the collection of maximal orders
(vertices) $\gamma \tilde \Lambda \gamma^{-1}$ have the same
geometric configuration as do the collection of $\tilde \Lambda \in
\C(\bnu)$, that is, they form a convex polytope in the apartment $\A = \gamma\A_0$.
\end{proof}
Finally, via propositions~\ref{prop:intersection-is-split} and
\ref{prop:split-is-intersection} we summarize the correspondence
between general split orders and convex polytopes in the building as
our main theorem.

\begin{thm}
  There is a one-to-one correspondence between convex polytopes (as described by
  Equation~\ref{eqn:linear-inequalities}) in apartments of the affine
  building for $SL_n(k)$  and split orders in $B = M_n(k)$.
\end{thm}


\bibliographystyle{amsplain} 

\providecommand{\bysame}{\leavevmode\hbox to3em{\hrulefill}\thinspace}
\providecommand{\MR}{\relax\ifhmode\unskip\space\fi MR }
\providecommand{\MRhref}[2]{%
  \href{http://www.ams.org/mathscinet-getitem?mr=#1}{#2}
}
\providecommand{\href}[2]{#2}

\end{document}